\documentclass[12pt]{article}

\usepackage[a4paper,margin=1in]{geometry}
\usepackage{amsmath,amssymb,amsfonts,amsthm}
\usepackage{graphicx}
\usepackage{hyperref}
\usepackage{bm}
\usepackage{enumerate}
\usepackage{mathrsfs}

\usepackage{tikz}
\usepackage{pgfplots}

\usetikzlibrary{arrows,shapes,positioning,calc}
\usepackage{subcaption}

\theoremstyle{definition}

\theoremstyle{remark}

\newcommand{\be}{\begin{equation}}
\newcommand{\ee}{\end{equation}}

\begin{document}

\title{Necessary and Sufficient Criterion for Singular or Nonsingular of Diagonally Dominant Matrices
\thanks{This paper is the English translation of the authors original article
published in \textit{Science China: Mathematics},
Vol.~44, No.~11, 2014. DOI: 10.1360/NO12014-00010.}
}
\author{Jidong Jin
\thanks{Department of Computer Science and Technology,Capital University of Economics and Business, Beijing 100070, China.}}
\date{November, 2014}
\maketitle

\begin{abstract}
The problem of determining whether a diagonally dominant matrix is singular or nonsingular is a classical topic in matrix theory. This paper develops necessary and sufficient conditions for the singularity or nonsingularity of diagonally dominant matrices. Starting from Taussky's theorem, we establish a unified line of theory which reduces the general problem to the study of irreducible diagonally dominant matrices. A complete similarity and unitary similarity analysis is given for singular irreducible diagonally dominant matrices, leading to a necessary and sufficient condition expressed in terms of an angle equation system associated with the nonzero off-diagonal entries. Applications and motivations from network dynamical systems, affine multi-agent systems, and Kolmogorov differential equations are also discussed.
\end{abstract}

\section{Introduction}

\

The problem of determining the nonsingularity of a matrix is a central topic in the theory of diagonally dominant matrices.
Taussky$^{[1]}$
surveyed the research carried out from the late 1870s to the 1940s in this area, while
Varga$^{[2]}$
reviewed the developments up to the mid-1970s.
The notions of generalized diagonal dominance and of ${\rm M}$C and ${\rm H}$Cmatrices introduced by
Ostrowski$^{[3-5]}$
have had a profound impact on the development of this field.
Berman and Plemmons$^{[6]}$
summarized fifty necessary and sufficient conditions for a matrix to be a nonsingular ${\rm M}$Cmatrix.
The characterization of ${\rm H}$Cmatrices has been one of the most active research topics in numerical analysis in China.
In this direction, Gao Yiming$^{[7-11]}$
Feng Mingxian$^{[12,13]}$
Huang Tingzhu$^{[14,15]}$
and
Li Wen$^{[16,17]}$
have all made many pioneering contributions.

The nonsingularity of a diagonally dominant matrix is closely related to the properties of the directed graph associated with the matrix.
Tausskys theorem$^{[17,18]}$
states that if the associated directed graph of a matrix is strongly connected, then the matrix is nonsingular as soon as it has at least one strictly diagonally dominant row.
Shivakumar and Chew$^{[20]}$
further proved that if, in the associated directed graph, every vertex can reach the vertex corresponding to some strictly diagonally dominant row, then the matrix is nonsingular.

The present paper studies necessary and sufficient conditions for the singularity or nonsingularity of diagonally dominant matrices.
In the early years of this century, Kolotilina$^{[21]}$
obtained a unitary similarity result for singular irreducible diagonally dominant matrices, which in effect provides a necessary and sufficient condition for the singularity of irreducible diagonally dominant matrices.
In this work we approach the problem from another angle, starting from Tausskys theorem, and derive more detailed results.
Moreover, Jin Jidong and Zheng Yufan$^{[22-24]}$
reduced the singularCnonsingular problem for general (reducible) diagonally dominant matrices to the corresponding problem for irreducible diagonally dominant matrices by using Tausskys theorem.
In other words, necessary and sufficient conditions for the singularity or nonsingularity of diagonally dominant matrices can be organized along a single line of theory stemming from Tausskys theorem.
Following this line, the present paper proceeds step by step and is organized as follows:

(1) We reduce the singularCnonsingular problem for reducible diagonally dominant matrices to the corresponding problem for irreducible diagonally dominant matrices;

(2) We analyze similarity and unitary similarity for singular irreducible diagonally dominant matrices;

(3) Finally, we present the main results of the paper: necessary and sufficient conditions for the singularity or nonsingularity of irreducible diagonally dominant matrices, together with a corresponding decision (testing) procedure.

\section{Taussky Theorem and Its Corollaries}

\

Unless otherwise specified, throughout this paper we assume that the matrix under consideration is
\(A=[a_{ij}]\in\mathbb{C}^{n\times n}\).
Here \(\mathbb{C}\) denotes the field of complex numbers and \(|\cdot|\) denotes the modulus of a complex number.

\textbf{Definition 1}\label{def:1}
(1) A matrix \(A=[a_{ij}]\in\mathbb{C}^{n\times n}\) is called
\emph{row diagonally dominant} if
\be\label{eq02-01}
\begin{array}{*{20}c}
{\left|a_{ii}\right|\ge \sum\limits_{\substack{j=1 \\ j\ne i}}^n \left|a_{ij}\right|}
&{}& i=1,\dots,n,
\end{array}
\ee
If
\[
\left|a_{ii}\right|>\sum_{\substack{j=1\\ j\ne i}}^n |a_{ij}|,
\]
then the \(i\)-th row is said to be \emph{strictly diagonally dominant};
if
\[
\left|a_{ii}\right|=\sum_{\substack{j=1\\ j\ne i}}^n |a_{ij}|,
\]
then the \(i\)-th row is \emph{weakly diagonally dominant}.
If all rows are strictly diagonally dominant, then \(A\) is called a
\emph{strictly row diagonally dominant matrix};
if all rows are weakly diagonally dominant, then \(A\) is called a
\emph{weakly row diagonally dominant matrix}.

(2) A matrix \(A=[a_{ij}]\in\mathbb{C}^{n\times n}\) is called
\emph{row balanced} if
\be\label{eq02-02}
\sum_{j=1}^n a_{ij}=0,\qquad i=1,\dots,n.
\ee
A row balanced matrix over the real field,
\(A=[a_{ij}]\in\mathbb{R}^{n\times n}\),
is called an \emph{affine matrix}.
Column diagonal dominance and column balance are defined analogously.

If a matrix \(A=[a_{ij}]\in\mathbb{C}^{n\times n}\) is diagonally dominant both by rows and by columns, it is called a
\emph{doubly diagonally dominant matrix};
if it is balanced both by rows and by columns, it is called a
\emph{doubly balanced matrix}.

\textbf{Definition 2}\label{def:2}
Given a matrix \(A=[a_{ij}]\in\mathbb{C}^{n\times n}\),
the simple directed graph
\(\Gamma(A)=\langle V, \mathcal{E}, A\rangle\)
is called the \emph{associated directed graph} of \(A\),
where \(V=\{v_1,\dots,v_n\}\) is the set of vertices
and
\[
\mathcal{E}
=\{(v_i,v_j)\in V\times V \mid a_{ij}\neq 0\}
\]
is the edge set.

\textbf{Definition 3}\label{def:3}
The associated directed graph \(\Gamma(A)=\langle V,\mathcal{E},A\rangle\)
is said to be \emph{strongly connected} if every vertex is reachable from every other vertex;
a matrix \(A\) is called \emph{irreducible} if its associated directed graph is strongly connected, and otherwise it is called \emph{reducible}.

The following classical result is due to Taussky$^{[18,19]}$.

\textbf{Theorem 2.1}\label{thm:2.1} (Taussky)
Let \(A=[a_{ij}]\in\mathbb{C}^{n\times n}\) be an irreducible diagonally dominant matrix.
If \(A\) has at least one strictly diagonally dominant row, then \(A\) is nonsingular.

\subsection{Frobenius Normal Form of Reducible Diagonally Dominant Matrices}

\

For any square matrix \(A\), there exists a permutation matrix \(P\) such that
\be\label{eq02-03}
P^\mathcal{T} A P
=
\left(
\begin{array}{*{20}c}
A_{11} &     &     &     &     &     \\
       &\ddots&    &     & 0   &     \\
0      &     &A_{ss}&    &     &     \\
A_{s+1,1}&\cdots&A_{s+1,s}&A_{s+1,s+1}&     &     \\
\vdots &\ddots&\vdots&\vdots&\ddots&     \\
A_{s+k,1}&\cdots&A_{s+k,s}&A_{s+k,s+1}&\cdots&A_{s+k,s+k}
\end{array}
\right).
\ee
The block matrix \((\ref{eq02-03})\) is called the \emph{Frobenius normal form} of \(A\).
The diagonal blocks \(A_{pp}\ (p=1,\dots,s+k)\) are irreducible square matrices
(the strongly connected components of \(\Gamma(A)\))
and are called the \emph{Frobenius blocks} of \(A\).
The Frobenius form \((\ref{eq02-03})\) has the following properties:
\be\label{eq02-04}
\left\{
\begin{array}{ll}
A_{pq}=0,&\forall\, q\ne p,\quad p\le s,\\[4pt]
A_{pq}\ne 0,&\exists\, q<p,\quad p>s.
\end{array}
\right.
\ee
Thus \(A_{11},\dots,A_{ss}\) are called the \emph{independent Frobenius blocks},
and
\(A_{s+1,s+1},\dots,A_{s+k,s+k}\) are called the \emph{non-independent Frobenius blocks}.

Jin and Zheng$^{[22-24]}$
obtained the following result.

\textbf{Corollary 2.1}\label{cor:2.1}
Every non-independent Frobenius block of a diagonally dominant matrix \(A\) is nonsingular.
The matrix \(A\) is singular if and only if at least one of its independent Frobenius blocks is singular.

\textbf{Proof.}

Suppose \((\ref{eq02-03})\) is the Frobenius normal form of \(A\)
and that \(A_{pp}\) is a non-independent Frobenius block.
Since \(A\) is diagonally dominant, so is \(A_{pp}\).
Because \(A_{pp}\) is a non-independent block, by \((\ref{eq02-04})\)
there exists \(q<p\) such that \(A_{pq}\neq 0\).
Therefore even if all rows of \(A_{pp}\) corresponding to the rows of \(A\) happen to be weakly diagonally dominant,
the block \(A_{pp}\) must contain at least one strictly diagonally dominant row.
As a Frobenius block, \(A_{pp}\) is irreducible, and hence by Tausskys theorem~\ref{thm:2.1} it is nonsingular.
Since all non-independent Frobenius blocks of \(A\) are nonsingular,
their corresponding rows in \(A\) are linearly independent.
Thus the singularity or nonsingularity of \(A\) is determined entirely by its independent Frobenius blocks:
\(A\) is singular if and only if at least one independent Frobenius block is singular. \(\Box\)

Corollary~\ref{cor:2.1} reduces the singularCnonsingular problem for general (reducible)
diagonally dominant matrices to the same problem for their independent Frobenius blocks.
Since Frobenius blocks are irreducible, this further reduces the problem to
determining the singularity or nonsingularity of irreducible diagonally dominant matrices.

\subsection{Irreducible Diagonally Dominant Matrices}

\

The following corollaries follow immediately from the definitions and from Taussky theorem~\ref{thm:2.1}.

\textbf{Corollary 2.2}\label{cor:2.2}
If an irreducible diagonally dominant matrix \(A\) is singular, then \(A\) must be weakly diagonally dominant.

\textbf{Corollary 2.3}\label{cor:2.3}
If an irreducible diagonally dominant matrix \(A=[a_{ij}]\in\mathbb{C}^{n\times n}\) is singular, then
\({\rm Rank}\,A = n-1\).

\textbf{Proof.}

Modify \(a_{ii}\) to some \(a'_{ii}\) such that \(|a'_{ii}| > |a_{ii}|\).
Then the modified matrix \(A'\) has its \(i\)-th row strictly diagonally dominant.
By Tausskys theorem~\ref{thm:2.1}, \(A'\) is nonsingular.
Therefore all rows of \(A\) except the \(i\)-th row must already be linearly independent.
Hence \({\rm Rank}\,A \ge n-1\).
Since \(A\) is singular, it follows that \({\rm Rank}\,A = n-1\).  \(\Box\)

This implies that for a singular irreducible diagonally dominant matrix \(A\in\mathbb{C}^{n\times n}\),
each row is a linear combination of the other \(n-1\) rows, and none of the coefficients in this combination can be zero.
We therefore obtain the following result.

\textbf{Corollary 2.4}\label{cor:2.4}
Let \(A\in\mathbb{C}^{n\times n}\) be a singular irreducible diagonally dominant matrix.
Let
\[
\rho=(\rho_1,\dots,\rho_n),\qquad
\gamma=(\gamma_1,\dots,\gamma_n)^\mathcal{T}
\]
be the left and right eigenvectors corresponding to the zero eigenvalue of \(A\), respectively.
Then neither \(\rho\) nor \(\gamma\) contains a zero entry.

\textbf{Definition 4}\label{def:4}
Let
\[
\rho=(\rho_1,\dots,\rho_n),\qquad
\gamma=(\gamma_1,\dots,\gamma_n)^\mathcal{T}
\]
be the left and right eigenvectors corresponding to the zero eigenvalue of
\(A\in\mathbb{C}^{n\times n}\).
The diagonal matrices
\[
P=\mathrm{diag}(\rho_1,\dots,\rho_n),\qquad
\Upsilon=\mathrm{diag}(\gamma_1,\dots,\gamma_n)
\]
are called, respectively, the \emph{left} and \emph{right eigenvector matrices} associated with the zero eigenvalue of \(A\).

By Corollary~\ref{cor:2.4}, we immediately obtain the following.

\textbf{Corollary 2.5}\label{cor:2.5}
For a singular irreducible diagonally dominant matrix \(A\), the left and right eigenvector matrices
\(P\) and \(\Upsilon\) associated with its zero eigenvalue are both nonsingular.

\section{Similarity Analysis of Singular Irreducible Diagonally Dominant Matrices}

\subsection{Balanced Diagonally Dominant Matrices}

\

The entries of a matrix \(A=[a_{ij}]\in\mathbb{C}^{n\times n}\) can be written in the polar form
\[
a_{ij} = \left|a_{ij}\right| e^{i\theta_{ij}},
\]
where \(\left|a_{ij}\right|\) is the modulus of \(a_{ij}\), \(e^{i\theta_{ij}}\) is the unit complex number associated with \(a_{ij}\), and \(\theta_{ij}\) is the argument of \(a_{ij}\).

\textbf{Proposition 3.1}\label{prop:3.1}
Let \(A=[a_{ij}]\in\mathbb{C}^{n\times n}\) be row diagonally dominant and row balanced.
Then

(1) \(A\) is weakly row diagonally dominant; and

(2) for every \(a_{ij}\neq 0\) with \(j\neq i\), one has
\(\theta_{ij} - \theta_{ii} = 2k\pi + \pi\) for some integer \(k\).

\textbf{Proof.}

Since \(A\) is row balanced, \eqref{eq02-02} holds.
Thus
\be\label{prop03-01-1}
- a_{ii} = \sum_{j=1,\; j\ne i}^n a_{ij}.
\ee
It follows that
\be\label{prop03-01-1-1}
\left|a_{ii}\right|
= \left|-a_{ii}\right|
= \left|\sum_{j=1,\; j\ne i}^n a_{ij}\right|
\le \sum_{j=1,\; j\ne i}^n \left|a_{ij}\right|.
\ee
Equality holds if and only if all nonzero \(a_{ij}\) \((a_{ij}\neq 0,\, j\neq i)\) have the same direction as vectors in the complex plane.

On the other hand, since \(A\) is row diagonally dominant,
\[
\left|a_{ii}\right|\ge \sum_{j=1,\; j\ne i}^n \left|a_{ij}\right|.
\]
Combining this with \eqref{prop03-01-1-1}, we obtain
\[
\left|a_{ii}\right|
= \sum_{j=1,\; j\ne i}^n \left|a_{ij}\right|.
\]
This proves (1).
Because equality in \eqref{prop03-01-1-1} holds, all arguments \(\theta_{ij}\) of the nonzero \(a_{ij}\) with \(j\neq i\) are the same.
Using \eqref{prop03-01-1}, it follows that
\(-a_{ii}\) must point in exactly the same direction as the sum of those \(a_{ij}\), hence the direction of \(a_{ii}\) differs from that of each such \(a_{ij}\) by \(\pi\) modulo \(2\pi\).
Thus
\(\theta_{ij} - \theta_{ii} = 2k\pi + \pi\) for some integer \(k\), which proves (2). \(\Box\)

\textbf{Proposition 3.2}\label{prop:3.2}
Let \(A=[a_{ij}]\in\mathbb{C}^{n\times n}\) be row balanced and column diagonally dominant.
Then

(1) \(A\) is weakly row diagonally dominant; and

(2) for every \(a_{ij}\neq 0\) with \(j\neq i\), one has
\(\theta_{ij} - \theta_{ii} = 2k\pi + \pi\) for some integer \(k\).

\textbf{Proof.}

Since \(A\) is row balanced, \eqref{eq02-02} holds and hence
\[
- a_{ii} = \sum_{j=1,\; j\ne i}^n a_{ij}.
\]
Therefore
\be\label{prop03-02-1}
\left|a_{ii}\right|
= \left|-a_{ii}\right|
= \left|\sum_{j=1,\; j\ne i}^n a_{ij}\right|
\le \sum_{j=1,\; j\ne i}^n \left|a_{ij}\right|.
\ee
Summing over \(i\) yields
\be\label{prop03-02-2}
\sum_{i=1}^n \left|a_{ii}\right|
\le \sum_{i=1}^n \sum_{j=1,\; j\ne i}^n \left|a_{ij}\right|.
\ee
Since \(A\) is column diagonally dominant, we have
\[
\left|a_{jj}\right|\ge \sum_{i=1,\; i\ne j}^n \left|a_{ij}\right|,
\]
and hence
\be\label{prop03-02-3}
\sum_{j=1}^n \left|a_{jj}\right|
\ge \sum_{j=1}^n \sum_{i=1,\; i\ne j}^n \left|a_{ij}\right|.
\ee
Comparing \eqref{prop03-02-2} and \eqref{prop03-02-3}, we see that equality must hold in \eqref{prop03-02-2}.
Moreover, equality in \eqref{prop03-02-2} can only occur when equality holds in \eqref{prop03-02-1} for every \(i\).
Thus \(A\) is weakly row diagonally dominant, which proves (1).

Since \(A\) is also row balanced, Proposition~\ref{prop:3.1} applies to give (2): for every \(a_{ij}\neq 0\) with \(j\neq i\), we have \(\theta_{ij} - \theta_{ii}=2k\pi + \pi\). \(\Box\)

By Proposition~\ref{prop:3.2}, we obtain the following result.

\textbf{Proposition 3.3}\label{prop:3.3}
Let \(A=[a_{ij}]\in\mathbb{C}^{n\times n}\) satisfy:
(1) \(A\) is diagonally dominant either by rows or by columns; and
(2) \(A\) is doubly balanced (both row balanced and column balanced).
Then

(1) \(A\) is weakly diagonally dominant both by rows and by columns; and

(2) for every \(a_{ij}\neq 0\) with \(j\neq i\),
\(\theta_{ij} - \theta_{ii} = 2k\pi + \pi\) and
\(\theta_{ij} - \theta_{jj} = 2k\pi + \pi\),
and hence \(\theta_{ii} - \theta_{jj} = 2k\pi\).

\textbf{Lemma 3.1}\label{lem:3.1}
Suppose a matrix \(A\) satisfies
(1) \(A\) is irreducible;
(2) \(A\) is diagonally dominant either by rows or by columns; and
(3) \(A\) is doubly balanced (row and column balanced).
Then all diagonal entries of \(A\) have the same argument, and each nonzero off-diagonal entry has an argument that differs from the common argument of the diagonal entries by \(2k\pi + \pi\).

\textbf{Proof.}

Since \(A\) is irreducible, its associated directed graph
\(\Gamma(A)=\langle V,\mathcal{E},A\rangle\) is strongly connected.
Hence there exists a set of edges
containing all vertices in a connected fashion:
\be\label{lem03-01-1}
\begin{array}{*{20}c}
\{(v_{s_1},v_{s_2}), (v_{s_2},v_{s_3}), \dots, (v_{s_{m-1}},v_{s_m})\} \subset \mathcal{E},
&,&
\{v_{s_1},v_{s_2},v_{s_3},\dots,v_{s_{m-1}},v_{s_m}\} = V.
\end{array}
\ee
The matrix \(A\) is diagonally dominant either by rows or by columns and is doubly balanced.
By Proposition~\ref{prop:3.3}, whenever \(a_{ij}\neq 0\), the diagonal entries \(a_{ii}\) and \(a_{jj}\) have the same argument.
By the definition of the associated directed graph, \(a_{ij}\neq 0\) is equivalent to \((v_i,v_j)\in\mathcal{E}\).
Combining this with \eqref{lem03-01-1}, we obtain
\(\theta_{s_1 s_1} = \cdots = \theta_{s_m s_m}\),
and this covers all diagonal entries of \(A\).
Applying Proposition~\ref{prop:3.3} again, we conclude that every nonzero off-diagonal entry of \(A\) has an argument that differs from the common argument of the diagonal entries by \(2k\pi + \pi\).
\(\Box\)

\subsection{Eigenvector Matrix Transformations}

\

A complex matrix \(A = [a_{ij}] \in \mathbb{C}^{n \times n}\) can be written as the Hadamard product of its modulus matrix
\(A_M = [|a_{ij}|] \in \mathbb{R}^{n \times n}\)
and its unitCcomplex matrix
\(A_C = [e^{i\theta(a_{ij})}] \in \mathbb{C}_0^{n \times n}\):
\[
A = [a_{ij}]
= \bigl[\,|a_{ij}| e^{i\theta(a_{ij})}\bigr]
= [A_M \circ A_C],
\]
where \(\mathbb{C}_0\) denotes the set of complex numbers of modulus one,
\(|a_{ij}|\) is the modulus of \(a_{ij}\), \(\theta(a_{ij})\) is its argument, and
\(e^{i\theta(a_{ij})}\) is the unitCcomplex factor of \(a_{ij}\).

\textbf{Definition 5}\label{def:5}
Notation convention: for a matrix \(A\), we write \(A_M\) for its modulus matrix and \(A_C\) for its unitCcomplex matrix.

\textbf{Lemma 3.2}\label{lem:3.2}
Let \(A = [a_{ij}]\) be a singular irreducible row diagonally dominant matrix.
Then there exist the left and right eigenvector matrices \(P\) and \(\Upsilon\) corresponding to the zero eigenvalue of \(A\) such that
\be\label{lem03-02-1}
P A \Upsilon = B,
\ee
where \(B \in \mathbb{R}^{n \times n}\) is a doubly balanced, doubly weakly diagonally dominant matrix.
Writing \(B = [B_M \circ B_C]\) with \(B_C = [b'_{ij}]\), we have
\be\label{lem03-02-2}
b'_{ij}
\begin{cases}
e^{i0} = 1, & j = i,\\[4pt]
e^{i\pi} = -1, & j \ne i,\ a_{ij} \ne 0,
\end{cases}
\ee
that is, the diagonal entries of \(B\) are positive (argument \(2k\pi\)), and the nonzero off-diagonal entries are negative (argument \(2k\pi + \pi\)).

\textbf{Proof.}

Since \(A\) is a singular irreducible row diagonally dominant matrix, Corollary~\ref{cor:2.5} of Tausskys theorem guarantees the existence of left and right eigenvector matrices \(P\) and \(\Upsilon\) corresponding to the zero eigenvalue of \(A\).
Both \(P\) and \(\Upsilon\) are nonsingular diagonal matrices.

Let \(A' = P A\).
Because \(P\) is diagonal, left-multiplication by \(P\) does not change the row diagonal dominance of \(A\), so \(A' = P A\) is again row diagonally dominant.
Moreover, since \(P\) is the left eigenvector matrix corresponding to the zero eigenvalue of \(A\), the matrix \(A' = P A\) is column balanced.
By Proposition~\ref{prop:3.2}, \(A' = P A\) is then also column diagonally dominant.
Thus \(A' = P A\) is column balanced and column diagonally dominant.

Next, let \(A'' = A' \Upsilon = P A \Upsilon\).
Because \(\Upsilon\) is the right eigenvector matrix corresponding to the zero eigenvalue of \(A\),
\(A'' = A' \Upsilon = P A \Upsilon\) is row balanced.
Furthermore, right-multiplication by \(\Upsilon\) does not affect the column balance or column diagonal dominance of \(A'\).
Hence \(A'' = P A \Upsilon\) is doubly balanced (both row and column balanced) and column diagonally dominant.

Since \(P\) and \(\Upsilon\) are diagonal with nonzero diagonal entries, the associated directed graphs \(\Gamma(A)\) and \(\Gamma(A'')\) are identical.
Because \(A\) is irreducible, so is \(A''\).
The matrix \(A'' = P A \Upsilon\) is doubly balanced, column diagonally dominant, and irreducible.
By Lemma~\ref{lem:3.1}, all diagonal entries of \(A''\) share the same argument, and every nonzero off-diagonal entry has an argument that differs from this common argument by \(2k\pi + \pi\).
Let \(\theta_D\) denote the common argument of the diagonal entries of \(A'' = P A \Upsilon\).
Then we can write
\[
P A \Upsilon = e^{i\theta_D} B,
\quad\text{so that}\quad
e^{-i\theta_D} P A \Upsilon = B.
\]
The matrix \(B \in \mathbb{R}^{n \times n}\) is therefore real.
Moreover, \(B\) inherits from \(A'' = P A \Upsilon\) the properties of being doubly balanced and doubly weakly diagonally dominant.
Since \(e^{-i\theta_D} P\) is again a left eigenvector matrix corresponding to the zero eigenvalue of \(A\), we may simply replace \(P\) by this scaled version and write
\[
P A \Upsilon = B.
\]
\(\Box\)

Lemma~\ref{lem:3.2} is the basis for the subsequent similarity analysis of singular irreducible diagonally dominant matrices.
It also has its own intrinsic significance, as it establishes a connection between singular irreducible diagonally dominant matrices and real, doubly weakly diagonally dominant matrices.

\subsection{Similarity}

\

A real matrix \(S = [s_{ij}] \in \mathbb{R}^{n \times n}\) is called a \emph{Markov matrix} if
\[
s_{ij} \ge 0 \quad (i,j = 1,\dots,n),
\qquad
\sum_{j=1}^n s_{ij} = 1 \quad (i=1,\dots,n).
\]
Let \(I\) denote the identity matrix.
Then \(I - S\) is a singular \(\mathbf{M}\)Cmatrix, called a \emph{Markov \(\mathbf{M}\)Cmatrix}.

\textbf{Definition 6}\label{def:6}
Notation convention:
for a square matrix \(A = [a_{ij}]_{n\times n}\), we define
\[
\mathcal{D}(A) = \mathrm{diag}(a_{11},\dots,a_{nn}),
\]
and call it the \emph{diagonal matrix} of \(A\).
We denote its inverse by \(\mathcal{D}^{-1}(A)\).

\textbf{Theorem 3.1}\label{thm:3.1}
Let \(A \in \mathbb{C}^{n\times n}\) be a singular irreducible row diagonally dominant matrix.
Then its right eigenvector matrix \(\Upsilon\) corresponding to the zero eigenvalue of \(A\) satisfies
\[
\Upsilon^{-1} A \Upsilon = \mathcal{D}(A)\,(I - S),
\]
where \(S\) is a Markov matrix.

\textbf{Proof.}

By Lemma~\ref{lem:3.2}, we have \(P A \Upsilon = B\).
On the left-hand side we perform the following equivalent transformation:
\[
P A \Upsilon
= P \bigl(\mathcal{D}(A)\mathcal{D}^{-1}(A)\bigr) A \Upsilon
= \bigl(P \mathcal{D}(A)\bigr)\bigl(\mathcal{D}^{-1}(A)A\bigr)\Upsilon.
\]
Thus
\[
\bigl(P \mathcal{D}(A)\bigr)\bigl(\mathcal{D}^{-1}(A)A\bigr)\Upsilon = B.
\]
Multiplying both sides on the left by \(\mathcal{D}^{-1}(B)\), we obtain
\be\label{thm03-001-001}
\bigl(\mathcal{D}^{-1}(B) P \mathcal{D}(A)\bigr)
\bigl(\mathcal{D}^{-1}(A)A\bigr)\Upsilon
= \mathcal{D}^{-1}(B) B.
\ee
The matrices \(\mathcal{D}^{-1}(A)A\) and \(\mathcal{D}^{-1}(B) B\) both have diagonal entries equal to \(1\).
Therefore
\[
\mathcal{D}^{-1}(B) P \mathcal{D}(A) = \Upsilon^{-1}.
\]
Substituting this expression into \eqref{thm03-001-001}, we get
\[
\Upsilon^{-1}\bigl(\mathcal{D}^{-1}(A)A\bigr)\Upsilon
= \mathcal{D}^{-1}(B) B.
\]
Since \(\Upsilon^{-1}\) and \(\mathcal{D}^{-1}(A)\) are diagonal, they commute:
\(\Upsilon^{-1}\mathcal{D}^{-1}(A) = \mathcal{D}^{-1}(A)\Upsilon^{-1}\).
Hence
\be\label{thm03-001-002}
\begin{array}{*{20}c}
\mathcal{D}^{-1}(A)\bigl(\Upsilon^{-1} A \Upsilon\bigr)
= \mathcal{D}^{-1}(B) B,
&;&
\Upsilon^{-1} A \Upsilon
= \mathcal{D}(A)\bigl(\mathcal{D}^{-1}(B) B\bigr).
\end{array}
\ee
Here \(B \in \mathbb{R}^{n\times n}\) is a real, doubly balanced matrix whose diagonal entries are positive and off-diagonal entries are negative.
Let \(\mathcal{D}^{-1}(B) B = [b'_{ij}]\).
Then
\[
b'_{ii} = 1,\qquad
-1 < b'_{ij} \le 0\ (j\ne i),
\qquad
\sum_{j=1,\; j\ne i}^n b'_{ij} = -1.
\]
Define a matrix \(S = [s_{ij}]\) by
\[
s_{ij} =
\begin{cases}
0, & j = i,\\[4pt]
- b'_{ij}, & j \ne i.
\end{cases}
\]
Then \(S\) is a Markov matrix, and
\(\mathcal{D}^{-1}(B) B = I - S\) is a Markov \(\mathbf{M}\)-matrix.
Substituting \(\mathcal{D}^{-1}(B) B = I - S\) into \eqref{thm03-001-002} yields
\[
\Upsilon^{-1} A \Upsilon = \mathcal{D}(A)(I - S).
\]
\(\Box\)

\subsection{Unitary similarity}

\

\textbf{Definition 7}\label{def:7}
For a diagonally dominant matrix \(A = [a_{ij}]\), its comparison matrix \(\mu(A) = [g_{ij}]\) is defined by
\[
g_{ij} =
\begin{cases}
|a_{ii}|, & j = i,\\[4pt]
-|a_{ij}|, & j \ne i.
\end{cases}
\]

By \eqref{lem03-02-2}, it is easy to see that
\be\label{Kolo2}
[A_M \circ B_C] = \mu(A).
\ee

In the following we will frequently use the following property of complex matrices:
let \(X, Y \in \mathbb{C}^{n \times n}\).
If at least one of \(X\) and \(Y\) is diagonal, then the product \(XY\) can be written as the Hadamard product of the products of their modulus matrices and unitCcomplex matrices, namely
\[
XY = [X_M \circ X_C][Y_M \circ Y_C]
= [X_M Y_M \circ X_C Y_C].
\]

\textbf{Theorem 3.2}\label{thm:3.2}
Let \(A \in \mathbb{C}^{n \times n}\) be a singular irreducible row diagonally dominant matrix.
Then the unitCcomplex matrix \(\Upsilon_C \in \mathbb{C}_0^{n \times n}\) (a diagonal unitary matrix) of the right eigenvector matrix \(\Upsilon\) corresponding to the zero eigenvalue of \(A\) satisfies
\be\label{thm03-02-01}
\Upsilon_C^{-1} A \Upsilon_C
= \mathcal{D}(A_C)[A_M \circ B_C]
= \mathcal{D}(A_C)\,\mu(A),
\ee
where \(B_C\) is the matrix satisfying \eqref{lem03-02-2}.

\textbf{Proof.}

By Lemma~\ref{lem:3.2} we have \(P A \Upsilon = B\).
Since \(P\) and \(\Upsilon\) are diagonal matrices, we may write
\[
\begin{array}{*{20}c}
P A \Upsilon
= [P_M \circ P_C][A_M \circ A_C][\Upsilon_M \circ \Upsilon_C]
= [P_M A_M \Upsilon_M \circ P_C A_C \Upsilon_C],
&;&
B = [B_M \circ B_C].
\end{array}
\]
Hence
\be\label{thm03-001-02}
P_C A_C \Upsilon_C = B_C.
\ee

Let \(\mathcal{D}(B_C)\) be the diagonal matrix formed from the diagonal entries of \(B_C\).
By \eqref{lem03-02-2}, we have \(\mathcal{D}(B_C) = I\).
Thus, taking the diagonal part in \eqref{thm03-001-02}, we obtain
\be\label{thm03-001-2}
\begin{array}{*{20}c}
P_C \mathcal{D}(A_C) \Upsilon_C = \mathcal{D}(B_C) = I,
&;&
P_C \mathcal{D}(A_C) = \Upsilon_C^{-1}.
\end{array}
\ee

Now consider \(P_C A \Upsilon_C\).
Since \(P_C\) and \(\Upsilon_C\) are diagonal matrices, we have
\[
P_C A \Upsilon_C
= [I \circ P_C][A_M \circ A_C][I \circ \Upsilon_C]
= [A_M \circ P_C A_C \Upsilon_C]
= [A_M \circ B_C].
\]
Thus
\[
P_C A \Upsilon_C = [A_M \circ B_C].
\]

Multiplying both sides of the last equality on the left by \(\mathcal{D}(A_C)\), and using the fact that
\(\mathcal{D}(A_C)\) and \(P_C\) are diagonal, we obtain
\[
\begin{array}{*{20}c}
\mathcal{D}(A_C)\bigl(P_C A \Upsilon_C\bigr)
= \mathcal{D}(A_C)[A_M \circ B_C],
&,&
\bigl(P_C \mathcal{D}(A_C)\bigr) A \Upsilon_C
= \mathcal{D}(A_C)[A_M \circ B_C].
\end{array}
\]
By \eqref{thm03-001-2}, we have \(P_C \mathcal{D}(A_C) = \Upsilon_C^{-1}\),
and therefore
\[
\Upsilon_C^{-1} A \Upsilon_C
= \mathcal{D}(A_C)[A_M \circ B_C].
\]
Finally, by \eqref{Kolo2} we have \([A_M \circ B_C] = \mu(A)\), so that
\[
\Upsilon_C^{-1} A \Upsilon_C
= \mathcal{D}(A_C)[A_M \circ B_C]
= \mathcal{D}(A_C)\,\mu(A).
\quad\Box
\]

\subsection{Unitary similarity of generalized diagonally dominant matrices}

\

In this subsection we discuss the relationship between Theorem~\ref{thm:3.2} and the result of Kolotilina$^{[21]}$.
For ease of comparison, we state the following corollary.

\textbf{Corollary 3.1}\label{cor:3.1}
Let \(A \in \mathbb{C}^{n \times n}\) be a singular irreducible row diagonally dominant matrix.
Then the unitCcomplex matrix \(\Upsilon_C \in \mathbb{C}_0^{n \times n}\) (a diagonal unitary matrix) of the right eigenvector matrix \(\Upsilon\) corresponding to the zero eigenvalue of \(A\) satisfies
\be\label{eq03-13}
\begin{array}{*{20}c}
\Upsilon_C^{-1} \mathcal{D}^{-1}(A) A \Upsilon_C
= \bigl[ \mathcal{D}^{-1}(A_M) A_M \circ B_C \bigr],
&,&
\Upsilon_C^{-1} \mathcal{D}^{-1}(A) A \Upsilon_C
= \mu\bigl( \mathcal{D}^{-1}(A) A \bigr).
\end{array}
\ee

\textbf{Proof.}

By Theorem~\ref{thm:3.2} we have
\[
\Upsilon_C^{-1} A \Upsilon_C = \mathcal{D}(A_C)\,[A_M \circ B_C].
\]
Multiplying both sides on the left by \(\mathcal{D}^{-1}(A)\), and using that
\(\mathcal{D}^{-1}(A)\), \(\Upsilon_C^{-1}\), and \(\mathcal{D}(A_C)\) are diagonal matrices, we obtain
\[
\mathcal{D}^{-1}(A)\Upsilon_C^{-1} A \Upsilon_C
= \Upsilon_C^{-1} \mathcal{D}^{-1}(A) A \Upsilon_C,
\]
and
\[
\mathcal{D}^{-1}(A)\mathcal{D}(A_C)\,[A_M \circ B_C]
= \bigl[ \mathcal{D}^{-1}(A_M) \circ \mathcal{D}^{-1}(A_C) \bigr]
   \bigl[ I \circ \mathcal{D}(A_C) \bigr]
   [A_M \circ B_C]
= \bigl[ \mathcal{D}^{-1}(A_M) A_M \circ B_C \bigr].
\]
Hence
\[
\Upsilon_C^{-1} \mathcal{D}^{-1}(A) A \Upsilon_C
= \bigl[ \mathcal{D}^{-1}(A_M) A_M \circ B_C \bigr].
\]
By \eqref{Kolo2} we have
\[
\bigl[ \mathcal{D}^{-1}(A_M) A_M \circ B_C \bigr]
= \mu\bigl( \mathcal{D}^{-1}(A) A \bigr).
\]
Substituting this into the previous equality yields
\[
\Upsilon_C^{-1} \mathcal{D}^{-1}(A) A \Upsilon_C
= \mu\bigl( \mathcal{D}^{-1}(A) A \bigr).
\quad\Box
\]

A square matrix \(A = [a_{ij}] \in \mathbb{C}^{n \times n}\) over the complex field is called a generalized diagonally dominant matrix
if there exist positive numbers \(v_1,\dots,v_n > 0\) such that
\be\label{eq1.2}
\begin{array}{*{20}c}
|a_{ii}|\, v_i \ge \displaystyle\sum_{j=1,\, j \ne i}^n |a_{ij}|\, v_j,
& & i = 1,\dots,n.
\end{array}
\ee

Kolotilina$^{[21]}$ proved a unitary similarity result for singular irreducible generalized diagonally dominant matrices (Theorem~2.3 in \cite{Kolotilina}).
We now derive this result from Corollary~\ref{cor:3.1}.

\textbf{Corollary 3.2} (Kolotilina)\label{cor:3.2}
Let \(A\) be a generalized diagonally dominant matrix satisfying \eqref{eq1.2}.
Then \(A\) is singular if and only if equality holds in \eqref{eq1.2} for every row of \(A\), and there exists a diagonal unitary matrix \(U\) such that
\be\label{Kolo}
U^{-1} \mathcal{D}^{-1}(A) A U = \mu\bigl( \mathcal{D}^{-1}(A) A \bigr).
\ee

\textbf{Proof.}

Since \(A\) is generalized diagonally dominant, there exists a nonsingular diagonal matrix
\[
V = \mathrm{diag}(v_1,\dots,v_n),
\quad v_1,\dots,v_n > 0,
\]
such that \(AV\) is diagonally dominant.
If \(A\) is singular, then \(AV\) is singular.
If \(A\) is irreducible, then \(AV\) is irreducible.
Substituting \(AV\) into \eqref{eq03-13} yields
\be\label{Kolo1}
\Upsilon_C^{-1} \bigl( \mathcal{D}^{-1}(AV)\,(AV) \bigr) \Upsilon_C
= \bigl[ \mathcal{D}^{-1}(A_M V)\,(A_M V) \circ B_C \bigr].
\ee

Since \(V\) is diagonal, we have \(\mathcal{D}(V) = V\).
For the left-hand side of \eqref{Kolo1}, we compute
\[
\Upsilon_C^{-1} \bigl( \mathcal{D}^{-1}(AV)\,(AV) \bigr) \Upsilon_C
= \Upsilon_C^{-1} \bigl( \mathcal{D}^{-1}(V)\,\mathcal{D}^{-1}(A)\,AV \bigr) \Upsilon_C
= \Upsilon_C^{-1} \bigl( V^{-1} \mathcal{D}^{-1}(A)\,AV \bigr) \Upsilon_C.
\]
Since \(\Upsilon_C, \Upsilon_C^{-1}, V^{-1}, V\) are all diagonal, we obtain
\[
\Upsilon_C^{-1} \bigl( V^{-1} \mathcal{D}^{-1}(A)\,AV \bigr) \Upsilon_C
= V^{-1} \Upsilon_C^{-1} \bigl( \mathcal{D}^{-1}(A)A \bigr) \Upsilon_C V.
\]
Multiplying on the left by \(V\) and on the right by \(V^{-1}\) gives
\[
\Upsilon_C^{-1} \bigl( \mathcal{D}^{-1}(A)A \bigr) \Upsilon_C.
\]

For the right-hand side of \eqref{Kolo1}, we have
\[
\bigl[ \mathcal{D}^{-1}(A_M V)\,(A_M V) \circ B_C \bigr]
= \bigl[ V^{-1} \mathcal{D}^{-1}(A_M) A_M V \circ B_C \bigr].
\]
Multiplying on the left by \(V\) and on the right by \(V^{-1}\), and using the diagonality of \(V\) and \(V^{-1}\), we obtain
\[
\begin{array}{*{20}l}
V \bigl[ V^{-1} \bigl( \mathcal{D}^{-1}(A_M) A_M \bigr) V \circ B_C \bigr] V^{-1}
&=& [V \circ I]\,[V^{-1} \bigl( \mathcal{D}^{-1}(A_M) A_M \bigr) V \circ B_C]\,[V^{-1} \circ I] \\[4pt]
&=& \bigl[ V V^{-1} \bigl( \mathcal{D}^{-1}(A_M) A_M \bigr) V V^{-1} \circ I B_C I \bigr] \\[4pt]
&=& \bigl[ \bigl( \mathcal{D}^{-1}(A_M) A_M \bigr) \circ B_C \bigr].
\end{array}
\]
Therefore,
\[
\Upsilon_C^{-1} \bigl( \mathcal{D}^{-1}(A)A \bigr) \Upsilon_C
= \bigl[ \bigl( \mathcal{D}^{-1}(A_M) A_M \bigr) \circ B_C \bigr].
\]
By \eqref{Kolo2},
\[
\bigl[ \bigl( \mathcal{D}^{-1}(A_M) A_M \bigr) \circ B_C \bigr]
= \mu\bigl( \mathcal{D}^{-1}(A) A \bigr),
\]
and hence
\[
\Upsilon_C^{-1} \bigl( \mathcal{D}^{-1}(A)A \bigr) \Upsilon_C
= \mu\bigl( \mathcal{D}^{-1}(A) A \bigr).
\]
Taking \(U = \Upsilon_C\) completes the proof. \(\Box\)

A diagonally dominant matrix is a special case of a generalized diagonally dominant matrix in the sense of Ostrowski, corresponding to the choice \(v = (1,\dots,1)^\mathcal{T}\) in \eqref{eq1.2}.
Thus Corollary~\ref{cor:3.1} is an immediate consequence of Corollary~\ref{cor:3.2}.
The above derivation, however, shows that the two are not merely in a one-way implication relation, but are in fact equivalent.
Corollary~\ref{cor:3.2} and Corollary~\ref{cor:3.1} form a pair of equivalent statements.

The equivalence between theorems in the theory of diagonally dominant matrices and those in the theory of generalized diagonally dominant matrices is not new.

\section{Necessary and sufficient conditions for the singularity and nonsingularity of irreducibly diagonally dominant matrices}

\

\textbf{Lemma 4.1}\label{lem:4.1}
Let $A \in \mathbb{C}^{n \times n}$ be an irreducible row diagonally dominant matrix.
If $A$ is singular, then the unitCcomplex matrix $\Upsilon_C \in \mathbb{C}_0^{n \times n}$ (a diagonal unitary matrix) associated with the right eigenvector matrix $\Upsilon$ corresponding to the zero eigenvalue of $A$ satisfies
\be\label{lem04-01-1}
\Upsilon_C^{-1} A_C \Upsilon_C = \mathcal{D}(A_C) B_C .
\ee
Here $A_C \in \mathbb{C}_0^{n \times n}$ is the unitCcomplex matrix of $A$, $\mathcal{D}(A_C)$ is the diagonal matrix formed by the diagonal entries of $A_C$, and $B_C$ is a matrix satisfying \textnormal{(\ref{lem03-02-2})}.

\textbf{Proof.}

By Theorem~\ref{thm:3.2}, we have
\[
\Upsilon_C^{-1} A \Upsilon_C = \mathcal{D}(A_C)[A_M \circ B_C] .
\]
On the one hand,
\[
\Upsilon_C^{-1} A \Upsilon_C
= \Upsilon_C^{-1}[A_M \circ A_C]\Upsilon_C
= [A_M \circ (\Upsilon_C^{-1} A_C \Upsilon_C)].
\]
On the other hand,
\[
\mathcal{D}(A_C)[A_M \circ B_C]
= [A_M \circ \mathcal{D}(A_C) B_C].
\]
Hence
\[
\Upsilon_C^{-1} A_C \Upsilon_C = \mathcal{D}(A_C) B_C .
\quad\Box
\]

We now state the main theorem of this paper: a necessary and sufficient condition for an irreducible diagonally dominant matrix to be singular.

\textbf{Theorem 4.1}\label{thm:4.1}
Let $A = [a_{ij}]$ be an irreducible diagonally dominant matrix.
Then $A$ is singular if and only if

(1) $A$ is weakly diagonally dominant; and

(2) there exist unit complex numbers
\[
\gamma_1 = e^{i\theta(\gamma_1)}, \dots, \gamma_n = e^{i\theta(\gamma_n)} \neq 0,
\]
such that the system of angle equations for the nonzero off-diagonal entries of $A$,
\be\label{thm04-01-1}
\bmod\bigl( -\theta(\gamma_i) + \theta(\gamma_j), 2\pi \bigr)
= \bmod\bigl( \pi + \theta(a_{ii}) - \theta(a_{ij}), 2\pi \bigr),
\qquad a_{ij} \neq 0,\ j \neq i,
\ee
or equivalently
\be\label{thm04-01-1a}
\bmod\bigl( \theta(\gamma_j), 2\pi \bigr)
= \bmod\bigl( \pi + \theta(a_{ii}) - \theta(a_{ij}) + \theta(\gamma_i), 2\pi \bigr),
\qquad a_{ij} \neq 0,\ j \neq i,
\ee
is consistent, i.e., $\theta(\gamma_1),\dots,\theta(\gamma_n)$ solve the above system.
Here $\theta(a_{ij})$ denotes the argument of $a_{ij}$, and $\theta(a_{ii})$ the argument of $a_{ii}$.
The notation $\bmod(x,y)$ denotes the remainder of $x$ upon division by $y$.

\textbf{Proof.}

(Necessity)
Suppose $A$ is an irreducible diagonally dominant matrix and $A$ is singular.
By Corollary~\ref{cor:2.2} of Tausskys theorem, $A$ is weakly diagonally dominant.

By Lemma~\ref{lem:4.1}, if $A$ is singular then
\be\label{thm04-01-03}
\Upsilon_C^{-1} A_C \Upsilon_C = \mathcal{D}(A_C) B_C .
\ee
Let $G = [g_{ij}] = \Upsilon_C^{-1} A_C \Upsilon_C$.
Then
\[
\begin{array}{*{20}l}
g_{ij} = \begin{cases}
e^{\,i(-\theta(\gamma_i) + \theta(a_{ii}) + \theta(\gamma_i))}, & j = i, \\
e^{\,i(-\theta(\gamma_i) + \theta(a_{ij}) + \theta(\gamma_j))}, & j \neq i,\ a_{ij} \neq 0,
\end{cases}
\\[8pt]
\theta(g_{ij}) =
\begin{cases}
\theta(a_{ii}), & j = i, \\
-\theta(\gamma_i) + \theta(a_{ij}) + \theta(\gamma_j), & j \neq i,\ a_{ij} \neq 0.
\end{cases}
\end{array}
\]

By Lemma~\ref{lem:4.1}, $B_C$ satisfies (\ref{lem03-02-2}).
Let $H = [h_{ij}] = \mathcal{D}(A_C) B_C$.
Then
\[
\begin{array}{*{20}l}
h_{ij} = \begin{cases}
e^{\,i(\theta(a_{ii}) + 2k\pi)}, & j = i, \\
e^{\,i(\theta(a_{ii}) + 2k\pi + \pi)}, & j \neq i,\ a_{ij} \neq 0,
\end{cases}
\\[8pt]
\theta(h_{ij}) =
\begin{cases}
2k\pi + \theta(a_{ii}), & j = i, \\
2k\pi + \pi + \theta(a_{ii}), & j \neq i,\ a_{ij} \neq 0.
\end{cases}
\end{array}
\]

From \eqref{thm04-01-03}, $G = H$, hence $g_{ij} = h_{ij}$ and thus $\theta(g_{ij}) = \theta(h_{ij})$.
Therefore, for any $j \neq i$ with $a_{ij} \neq 0$,
\[
\begin{array}{*{20}c}
-\theta(\gamma_i) + \theta(a_{ij}) + \theta(\gamma_j)
= 2k\pi + \pi + \theta(a_{ii}),
&;&
-\theta(\gamma_i) + \theta(\gamma_j)
= 2k\pi + \pi + \theta(a_{ii}) - \theta(a_{ij}).
\end{array}
\]
Eliminating the integer parameter $k$ via the modulo operation yields \eqref{thm04-01-1}.
Since
\[
\Upsilon_C = \mathrm{diag}(\gamma_1,\dots,\gamma_n)
= \mathrm{diag}\bigl(e^{i\theta(\gamma_1)},\dots,e^{i\theta(\gamma_n)}\bigr)
\]
is the unitCcomplex matrix of the right eigenvector matrix $\Upsilon$ corresponding to the zero eigenvalue of $A$, and by Corollary~\ref{cor:2.5} of Tausskys theorem $\Upsilon$ is nonsingular, we have
\[
\gamma_1 = e^{i\theta(\gamma_1)}, \dots, \gamma_n = e^{i\theta(\gamma_n)} \neq 0.
\]

(Sufficiency)
Conversely, assume there exist
\[
\gamma_1 = e^{i\theta(\gamma_1)}, \dots, \gamma_n = e^{i\theta(\gamma_n)} \neq 0
\]
such that \eqref{thm04-01-1} holds.
Set
\[
Q = \mathrm{diag}(\gamma_1,\dots,\gamma_n)
= \mathrm{diag}\bigl(e^{i\theta(\gamma_1)},\dots,e^{i\theta(\gamma_n)}\bigr).
\]
Then $Q$ is a diagonal unitary matrix.
Consider $Q^{-1} A Q$.
Note that $Q^{-1}$ and $\mathcal{D}(A_C)$ are diagonal matrices, and
\[
Q^{-1} A Q
= Q^{-1} \bigl( \mathcal{D}(A_C)\mathcal{D}^{-1}(A_C) \bigr) A Q
= \mathcal{D}(A_C)\bigl( Q^{-1} \mathcal{D}^{-1}(A_C) A Q \bigr).
\]

We now show that $Q^{-1} \mathcal{D}^{-1}(A_C) A Q$ is singular.
Since $Q$ is a unitCcomplex diagonal matrix, and using the diagonality of $\mathcal{D}^{-1}(A_C)$, $Q^{-1}$, and $Q$, we have
\[
Q^{-1} \mathcal{D}^{-1}(A_C) A Q
= [I \circ Q^{-1}] [I \circ \mathcal{D}^{-1}(A_C)] [A_M \circ A_C] [I \circ Q]
= [A_M \circ Q^{-1} \mathcal{D}^{-1}(A_C) A_C Q].
\]
Let $B'_C = Q^{-1} \mathcal{D}^{-1}(A_C) A_C Q = [b_{ij}]$.
We first show that $B'_C$ satisfies (\ref{lem03-02-2}).
For the diagonal entries,
\[
b_{ii}
= e^{\,i\bigl(-\theta(\gamma_i) - \theta(a_{ii}) + \theta(a_{ii}) + \theta(\gamma_i)\bigr)}
= e^{i0} = 1,
\]
so $b_{ii}$ satisfies (\ref{lem03-02-2}).
For $j \neq i$ with $a_{ij} \neq 0$,
\[
b_{ij}
= e^{\,i\bigl(-\theta(\gamma_i) - \theta(a_{ii}) + \theta(a_{ij}) + \theta(\gamma_j)\bigr)}
= e^{\,i\bigl(-\theta(\gamma_i) + \theta(\gamma_j) - \theta(a_{ii}) + \theta(a_{ij})\bigr)}.
\]
By \eqref{thm04-01-1},
\[
\bmod\bigl( -\theta(\gamma_i) + \theta(\gamma_j) - \theta(a_{ii}) + \theta(a_{ij}), 2\pi \bigr)
= \pi.
\]
Hence
\[
\bmod\bigl( \theta(b_{ij}), 2\pi \bigr)
= \bmod\bigl( -\theta(\gamma_i) + \theta(\gamma_j) - \theta(a_{ii}) + \theta(a_{ij}), 2\pi \bigr)
= \pi,
\]
which implies $b_{ij} = -1$, so $b_{ij}$ also satisfies (\ref{lem03-02-2}).

Since $A_M = [|a_{ij}|]$ is the modulus matrix of $A$ and $B'_C$ satisfies (\ref{lem03-02-2}), we have
\[
Q^{-1} \mathcal{D}^{-1}(A_C) A Q
= [A_M \circ B'_C] = [a'_{ij}],
\]
where
\[
a'_{ii} = |a_{ii}|, \qquad a'_{ij} = -|a_{ij}| \quad (j \neq i).
\]
Because $A$ is weakly diagonally dominant,
\[
|a_{ii}| = \sum_{j=1,\, j \ne i}^n |a_{ij}|,
\]
we obtain
\[
a'_{ii} = |a_{ii}| = \sum_{j=1,\, j \ne i}^n |a_{ij}|
= \sum_{j=1,\, j \ne i}^n \bigl(-a'_{ij}\bigr),
\qquad
\sum_{j=1}^n a'_{ij} = 0.
\]
Thus $Q^{-1} \mathcal{D}^{-1}(A_C) A Q$ is a rowCbalanced matrix and therefore singular.
Since $\mathcal{D}(A_C)$ is nonsingular (the diagonal entries of $A$ are nonzero for a diagonally dominant matrix), it follows that $A$ itself must be singular.
\(\Box\)

Regarding the singularity of irreducible diagonally dominant matrices, Tausskys theorem~\ref{thm:2.1} provides a necessary condition: if an irreducible diagonally dominant matrix is singular, then it must in fact be weakly diagonally dominant.
Theorem~\ref{thm:4.1} strengthens this by adding a condition on the existence (consistency) of a solution of the angle equation system \eqref{thm04-01-1} for the nonzero off-diagonal entries, thereby yielding a necessary and sufficient condition.
In this sense, Theorem~\ref{thm:4.1} refines and extends Tausskys theorem.

Given a matrix $A$, the arguments $\theta(a_{ii})$ and $\theta(a_{ij})$ are fixed.
Each equation in \eqref{thm04-01-1} is in fact a two-variable equation with respect to
$\theta(\gamma_i)$ and $\theta(\gamma_j)$.
Every nonzero off-diagonal entry $a_{ij}$ of $A$ corresponds to one such equation.
For a given irreducible weakly diagonally dominant matrix $A$, one may determine the
singularity or nonsingularity of $A$ by solving the system \eqref{thm04-01-1}.
Theorem~\ref{thm:4.1} shows that if \eqref{thm04-01-1} has a solution, then $A$ is singular;
otherwise, $A$ is nonsingular.

\textbf{Corollary 4.1}\label{cor:4.1}
Suppose
\[
\gamma_1 = e^{i\theta(\gamma_1)}, \dots,
\gamma_n = e^{i\theta(\gamma_n)} \neq 0
\]
is a solution of the system \eqref{thm04-01-1}.
Let $\beta$ be any unit complex number.
Then
\be\label{prop04-01-1}
\beta \gamma_1 = e^{i\bigl(\theta(\gamma_1) + \theta(\beta)\bigr)}, \dots,
\beta \gamma_n = e^{i\bigl(\theta(\gamma_n) + \theta(\beta)\bigr)}
\ee
is also a solution of the system \eqref{thm04-01-1}.

\textbf{Proof.}

Substituting \eqref{prop04-01-1} into \eqref{thm04-01-1} gives
\[
\bmod\Bigl(
\theta(a_{ij}) - \theta(a_{ii})
+ \bigl(\theta(\gamma_j) + \theta(\beta)\bigr)
- \bigl(\theta(\gamma_i) + \theta(\beta)\bigr),
2\pi
\Bigr)
= \pi .
\]
This is equivalent to
\[
\bmod\bigl(
\theta(a_{ij}) - \theta(a_{ii})
+ \theta(\gamma_j) - \theta(\gamma_i),
2\pi
\bigr)
= \pi,
\]
which is precisely \eqref{thm04-01-1}.
\(\Box\)

This corollary shows that solving the system \eqref{thm04-01-1} does not require
general-purpose methods for solving systems of equations; instead, one can determine
the variables successively by substitution.
For instance, one may fix a single angle, say
$\theta(\gamma_1) = 0$ (so that $\gamma_1 = e^{i0} = 1$), and then use
\eqref{thm04-01-1a} iteratively to compute
$\theta(\gamma_2), \dots, \theta(\gamma_n)$ one by one.

For a complex number $a = |a| e^{i\theta}$, if $\theta \in \{0,\pi\}$ then $a$ is real.
Thus we call
\[
\Pi = \{0,\pi\}
\]
the \emph{realCangle domain}.
The solving process in Section~4.2 shows that if $A$ is a real matrix, then when $\theta(\gamma_1)=0$, we have
$\theta(\gamma_1),\dots,\theta(\gamma_n)\in\{0,\pi\}$.
This is easy to see from the solving formula \eqref{thm04-01-1a}:
\be\label{real}
\theta(\gamma_j)
= \mathrm{mod}\!\left(\pi + \theta(a_{ii}) - \theta(a_{ij}) + \theta(\gamma_i),\,2\pi\right).
\ee
Since $A$ is a real matrix, $\theta(a_{ii}),\theta(a_{ij})\in\{0,\pi\}$.
Therefore, as long as $\theta(\gamma_i)\in\{0,\pi\}$, one must have $\theta(\gamma_j)\in\{0,\pi\}$.
Fixing any $\theta(\gamma_1)\in\{0,\pi\}$, the solving process ensures that
$\theta(\gamma_j)\in\{0,\pi\}$ for all $j=2,\dots,n$.
In fact, choosing any $\theta(\gamma_i)\in\{0,\pi\}$ $(i=1,\dots,n)$ as the initial value makes no difference.
Hence the following corollary follows immediately.

\textbf{Corollary 4.2}\label{cor:4.2}
If $A\in\mathbb{R}^{n\times n}$ is a singular irreducible (weakly) diagonally dominant matrix over the real field,
then the angle equation system \eqref{thm04-01-1} admits a solution in the realCangle domain; that is,
there exist
\be\label{cor04-02-0001}
\theta(\gamma_1),\dots,\theta(\gamma_n)\in\{0,\pi\}
\ee
which solve the angle system \eqref{thm04-01-1}.
Moreover, as long as there exists any index $i\in\{1,\dots,n\}$ such that $\theta(\gamma_i)\in\{0,\pi\}$,
the condition \eqref{cor04-02-0001} holds.

Since \eqref{cor04-02-0001} specifies the arguments of the entries of the unitCcomplex similarity matrix $\Upsilon_C$,
we obtain the following corollary.

\textbf{Corollary 4.3}\label{cor:4.3}
If $A\in\mathbb{R}^{n\times n}$ is a singular irreducible diagonally dominant matrix,
then there exists a real unitary similarity transformation matrix
\[
\Upsilon_C
= \operatorname{diag}(\gamma_1,\dots,\gamma_n)
= \operatorname{diag}\!\left(e^{i\theta(\gamma_1)},\dots,e^{i\theta(\gamma_n)}\right)
\in\mathbb{R}^{n\times n}
\]
which satisfies \eqref{thm03-02-01},
where $\theta(\gamma_1),\dots,\theta(\gamma_n)\in\{0,\pi\}$.

\section{Discussion and Conclusion}

\

One methodological feature of the present paper is that it starts from the
properties of singular diagonally dominant matrices.

This section explains that this choice is not merely methodological. The work
grew out of a concrete mathematical problem encountered by the author in the
study of network dynamics, namely the structure of singular diagonally
dominant matrices.

In fact, network dynamics has already led the theory of diagonally dominant
matrices towards a new and important direction: the study of diagonally
dominant matrices whose associated directed graph structure plays a crucial
role. In this section, we briefly discuss the classes of problems and
application domains that motivate this point of view.

\subsection{Problems and domains related to singular diagonally dominant matrices}

\subsubsection{Network dynamical systems in an inertial frame}

\

In recent years, multi-agent systems have attracted extensive attention from
researchers in systems science, control theory, intelligent control, computer
science, communication networks, biology, sociology, and related fields.
Among these problems, multi-body dynamics in a Galilean inertial frame forms
one of the most important classes.

Let $O$ be the set of agents and $o_i \in O$ a generic agent. For each
time~$t$, let $N(t,o_i)$ denote the set of all agents that \emph{act} on
$o_i$ at time $t$; this is called the neighbour set of $o_i$ at time $t$.
The family of neighbour sets of all agents induces a time-varying directed
interaction network $\mathbb{G}(t)$:
{\small
\be
\begin{array}{*{20}c}
   {\mathbb{G}(t) =  < O,{\cal E}(t) > } & {;} &
   {\mathcal{E}(t) = \left\{ {\left. {(o_j ,o_i )_t  \in O \times O} \right|
   o_i  \in O,\ o_j  \in N(t,o_i )} \right\}}  \\
\end{array} .
\ee}

Jin Jidong and Zheng Yufan$^{[25\text{--}30]}$

studied the following velocity coordination model for multi-body dynamical
systems in an inertial frame:
{\small
\be\label{eq010301}
{\dot v_i (t) = \sum\limits_{o_j  \in N(t,o_i )} { - g_{ij} \left( {\left\| {v_i (t) - v_j (t)} \right\|} \right)} \vec v_{ij}
 = \sum\limits_{o_j  \in N(t,o_i )} { - g_{ij} \left( {\left\| {v_i (t) - v_j (t)} \right\|} \right)\frac{{v_i (t) - v_j (t)}}{{\left\| {v_i (t) - v_j (t)} \right\|}}} }
\ee}
where $v_i(t)$ and $\dot v_i(t)$ denote the velocity and acceleration of
agent $o_i$ at time $t$, respectively; $g_{ij}(\bullet) \ge 0$ is a scalar
function representing the interaction strength from $o_j$ to $o_i$; and
$\vec{v}_{ij}=\dfrac{v_i (t) - v_j (t)}{\left\| {v_i (t) - v_j (t)} \right\|}$
is the unit vector in the direction of $v_i(t)-v_j(t)$.
The minus sign in front of $g_{ij}(\cdot)$ reflects the tendency of $v_i(t)$
to approach $v_j(t)$.

Model~\eqref{eq010301} is a general mechanical model describing multi-body
dynamics in an inertial frame. It differs from classical mechanical models in
that interactions between agents are not required to be reciprocal: the
influence from one agent to another may be purely one-directional. This
framework is therefore well suited for describing velocity coordination in
multi-body systems composed of natural or artificial autonomous agents
(such as bird flocks, fish schools, traffic flow, unmanned aerial vehicles
and robotic swarms) moving in Galilean space.

Define
\[
{g'_{ij} (t) = \frac{{g_{ij} \left( {\left\| {v_i (t) - v_j (t)} \right\|} \right)}}{{\left\| {v_i (t) - v_j (t)} \right\|}}}
\]
and the linear operator $\mathcal{Q}(t)=[\mathcal{G}_{ij}(t)]_{n \times n}$ by
{\small
\be\label{eq0106}
\mathcal{G}_{ij} (t) = \left\{ {\begin{array}{*{20}l}
   {  \sum\limits_{o_j  \in N(t,o_i )} {-g'_{ij} (t)} } & {j = i}  \\[0.3em]
   {g'_{ij} (t)} & {j \ne i}  \\
\end{array}} \right.
\ee}
Then~\eqref{eq010301} can be written compactly as
{\small
\[
\dot v(t) = \left( {\mathcal{Q}(t) \otimes I_m } \right)\left( {v(t) \otimes I_m } \right),
\]}
where $v(t) = \bigl( v_1 (t), \ldots ,v_n (t) \bigr)^\mathcal{T}$,
$\otimes$ denotes the Kronecker product, and $I_m$ $(m=1,2,3)$ is the
$m$-dimensional identity matrix.

From~\eqref{eq0106} it follows that $\mathcal{Q}(t)$ and
$\mathcal{Q}(t) \otimes I_m$ form a special class of singular diagonally
dominant matrices: they are row-balanced diagonally dominant matrices; and
since they are real matrices, they are affine diagonally dominant matrices.
Consequently, the multi-body dynamical model~\eqref{eq010301} in an inertial
frame is a first-order dynamical system governed by an affine diagonally
dominant linear operator.

\subsubsection{Linear network dynamical systems}

\

If there exist constants $\gamma'_{ij}$ such that, for all $y \ge 0$,
\[
\gamma '_{ij}  = \frac{{g_{ij} (y)}}{y}, \qquad i,j = 1, \ldots ,n, \ j \ne i,
\]
then~\eqref{eq010301} reduces to a first-order linear time-invariant system
{\small
\be\label{eq010901}
\begin{array}{*{20}l}
   {\dot x(t) = \Gamma x(t)} & {} &
   { \gamma _{ij}  = \left\{ {\begin{array}{*{20}c}
   {\displaystyle\sum\limits_{o_j  \in N(t,o_i )} {-\gamma '_{ij} } } & {j = i}  \\[0.4em]
   { \gamma '_{ij} } & {j \ne i}  \\
\end{array}} \right.}  \\
\end{array}
\ee}
From~\eqref{eq010301} we see that $\gamma'_{ij} \ge 0$ for all $i,j$
with $j \neq i$, hence $\Gamma=[\gamma_{ij}]$ is an affine diagonally
dominant matrix.

The study of the non-autonomous/autonomous discrete-time versions of
system~\eqref{eq010901}, namely
{\small
\[
x(t + 1) = S(t)x(t),
\]}
where the linear operator $S(t)$ is a Markov matrix with nonzero diagonal
entries, dates back to the work of DeGroot$^{[31]}$

on consensus formation in collective decision making, and the work of
Tsitsiklis and Bertsekas$^{[32,33]}$

on cooperative optimal control in multiprocessor systems. The former concerns
problems in sociology and economics; the latter provides theoretical
guidance for contemporary information science topics, such as optimal
control of intelligent communication networks and cloud (computing/storage)
systems.

\subsubsection{Kolmogorov differential equations}

\

In fact, singular diagonally dominant matrices arise not only in network
dynamics.

Let $Q=[q_{ij}]$ be a Markov transition matrix,
\[
0 \le q_{ij} \le 1, \qquad \sum_{j = 1}^n q_{ij} = 1,
\]
and let $P(t)=[p_{ij}(t)]$ denote the transition probability matrix of a
continuous-time Markov process at time $t$. Define $\Gamma=Q-I$. Then
{\small
\be\label{eq01010}
\dot P(t) = \Gamma P(t) = (Q - I)P(t)
\ee}
is the classical Kolmogorov forward differential equation. Here
$\Gamma=Q-I$ is an affine diagonally dominant matrix, and the equation
\eqref{eq01010} can be viewed as a special case of the first-order linear
time-invariant system~\eqref{eq010901}.

\subsection{Concluding remarks}

\

Starting from Taussky theorem and following the line of inquiry dictated by
the properties of singular diagonally dominant matrices, this paper derives a
necessary and sufficient condition for the singularity or nonsingularity of
diagonally dominant matrices. As emphasized above, this is not merely a
matter of methodology. Singular diagonally dominant matrices constitute a
basic mathematical form underlying a large class of closely related problems
in mathematics, the natural sciences, and modern engineering.

A particular emphasis of this work is the role played by the associated
directed graph of a matrix in the analysis of diagonally dominant matrices.

The result of this paper stating that the singularity or nonsingularity of a
reducible diagonally dominant matrix is determined solely by its independent
Frobenius blocks (Corollary~\ref{cor:2.1} of Taussky's theorem) provides concrete
support for the following viewpoint: without understanding the structure of
the associated directed graph of a matrix, some questions in matrix theory
cannot even be formulated clearly; once this structure is properly taken into
account, certain long-standing conceptual difficulties may become almost
trivial.


\bigskip







\begin{thebibliography}{99}

\bibitem{Taussky1}
Taussky O.
A recurring theorem on determinants.
\emph{The American Mathematical Monthly}, 1949, \textbf{56}(10): 672--676.

\bibitem{Varga}
Varga R.~S.
On recurring theorems on diagonal dominance.
\emph{Linear Algebra and Its Applications}, 1976, \textbf{13}(1): 1--9.

\bibitem{Ostrowski1}
Ostrowski A.~M.
\"Uber die Determinanten mit \"uberwiegender Hauptdiagonale.
\emph{Commentarii Mathematici Helvetici}, 1937, \textbf{10}(1): 69--96.

\bibitem{Ostrowski2}
Ostrowski A.~M.
Determinanten mit \"uberwiegender Hauptdiagonale und die absolute Konvergenz von linearen Iterationsprozessen.
\emph{Commentarii Mathematici Helvetici}, 1956, \textbf{30}(1): 175--210.

\bibitem{Ostrowski3}
Ostrowski A.~M.
On some metrical properties of operator matrices and matrices partitioned into blocks.
\emph{Journal of Mathematical Analysis and Applications}, 1961, \textbf{2}(2): 161--209.

\bibitem{Berman}
Berman A., Plemmons R.~J.
\emph{Nonnegative Matrices in the Mathematical Sciences}.
New York: Academic Press, 1979.

\bibitem{Gao1}
Gao Y.~M.
Criteria for generalized diagonally dominant matrices and nonsingular matrices.
\emph{Journal of Northeast Normal University}, 1982, \textbf{3}: 23--28. (in Chinese)

\bibitem{Gao2}
Gao Y.~M.
Criteria for generalized diagonal dominance and nonsingularity of matrices (II).
\emph{Journal of Engineering Mathematics}, 1988, \textbf{5}(3): 12--16. (in Chinese)

\bibitem{Gao3}
Gao Y.~M.
Criteria for generalized diagonally dominant matrices and $M$-matrices.
\emph{Journal of Computational Mathematics in Colleges and Universities}, 1992, \textbf{14}(3): 233--239. (in Chinese)

\bibitem{Gao4}
Gao Y.~M., Wang X.~H.
Criteria for generalized diagonally dominant matrices and $M$-matrices.
\emph{Linear Algebra and Its Applications}, 1992, \textbf{169}: 257--268.

\bibitem{Gao5}
Gao Y.~M., Wang X.~H.
Criteria for generalized diagonally dominant matrices and $M$-matrices. II.
\emph{Linear Algebra and Its Applications}, 1996, \textbf{248}: 339--353.

\bibitem{feng1}
Feng M.~X.
Criteria for generalized diagonally dominant matrices and their applications.
\emph{Mathematica Annals}, 1985, \textbf{6A}(3): 323--330. (in Chinese)

\bibitem{feng2}
Feng M.~X.
Generalizations of diagonal dominance of matrices and their applications.
\emph{Acta Mathematicae Applicatae Sinica}, 1989, \textbf{12}(1): 35--43. (in Chinese)

\bibitem{Huang2}
Huang T.~Z.
Simple criteria for nonsingular $H$-matrices.
\emph{Journal of Computational Mathematics}, 1993, \textbf{15}(3): 318--328. (in Chinese)

\bibitem{Huang3}
Huang T.~Z.
A generalization of Ostrowski's theorem and conditions for nonsingular $H$-matrices.
\emph{Journal of Computational Mathematics}, 1994, \textbf{16}(1): 19--24. (in Chinese)

\bibitem{Li1}
Li W.
Criteria for generalized diagonally dominant matrices.
\emph{Journal of Applied Mathematics and Computation}, 1995, \textbf{9}(2): 35--38. (in Chinese)

\bibitem{Li2}
Li W.
A note on criteria for generalized diagonally dominant matrices.
\emph{Journal of Computational Mathematics in Colleges and Universities}, 1997, \textbf{1}: 93--96. (in Chinese)

\bibitem{Taussky2}
Taussky O.
Bounds for characteristic roots of matrices.
\emph{Duke Mathematical Journal}, 1948, \textbf{15}(4): 1043--1044.

\bibitem{Taussky3}
Taussky O.
Bounds for the characteristic roots of matrices II.
\emph{Journal of Research of the National Bureau of Standards}, 1951, \textbf{46}(2): 124--125.

\bibitem{Shivakumar}
Shivakumar P.~N., Chew K.~H.
A sufficient condition for nonvanishing of determinants.
\emph{Proceedings of the American Mathematical Society}, 1974, \textbf{44}(1): 63--66.

\bibitem{Kolotilina}
Kolotilina L.~Y.
Nonsingularity/singularity criteria for nonstrictly block diagonally dominant matrices.
\emph{Linear Algebra and Its Applications}, 2003, \textbf{359}(1--3): 133--159.

\bibitem{Jin1}
Jin J.~D., Zheng Y.~F.
Consensus of multi-agent system under directed network:
A matrix analysis approach.
In: \emph{Proceedings of the IEEE International Conference on Control and Automation (ICCA 2009)}.
Christchurch: IEEE, 2009, 280--284.

\bibitem{Jin2}
Jin J.~D., Zheng Y.~F.
Global behavior of multi-agent systems under directed networks.
\emph{Systems Science and Mathematics}, 2011, \textbf{31}(1): 114--122. (in Chinese)

\bibitem{Jin3}
Jin J.~D.
Collective behavior of cooperative multi-agent dynamical systems.
Ph.D. Thesis, Shanghai University, 2011. (in Chinese)

\bibitem{Jin4}
Jin J.~D., Zheng Y.~F.
The collective behavior of asymmetric affine multi-agent system.
In: \emph{Proceedings of the 8th Asian Control Conference (ASCC 2011)}.
Kaohsiung: IEEE, 2011, 800--805.

\bibitem{Jin5}
Jin J.~D., Zheng Y.~F.
Global behavior of affine nonlinear multi-agent cooperative dynamical systems under directed networks.
\emph{Control Theory and Applications}, 2011, \textbf{28}(10): 1377--1383. (in Chinese)

\bibitem{Jin6}
Jin J.~D., Zheng Y.~F., Zheng X.~L.
A unified theory for collective behavior of cooperative systems.
In: \emph{Proceedings of the 9th IEEE International Conference on Control and Automation (ICCA 2011)}.
Santiago: IEEE, 2011, 471--476.

\bibitem{Jin6-1}
Zheng Y.~F., Jin J.~D.
Necessary condition of consensus for affine multi-agent systems under time-varying directed networks.
In: \emph{Proceedings of the 25th Chinese Control and Decision Conference (CCDC 2013)}.
Guiyang: IEEE, 2013, 702--706.

\bibitem{Jin6-2}
Jin J.~D., Zheng Y.~F.
Necessary and sufficient condition of consensus for affine multi-agent cooperative systems under time-varying directed networks.
In: \emph{Proceedings of the 9th Asian Control Conference (ASCC 2013)}.
Istanbul: IEEE, 2013, 1--6.

\bibitem{Jin7}
Jin J.~D., Zheng Y.~F.
State consensus of affine multi-agent systems under time-varying directed networks.
\emph{Science China Information Sciences}, 2013, \textbf{43}(10): 1365--1382. (in Chinese)

\bibitem{DeGroot}
DeGroot M.~H.
Reaching a consensus.
\emph{Journal of the American Statistical Association}, 1974, \textbf{69}(345): 118--121.

\bibitem{Tsitsiklis1}
Tsitsiklis J.~N.
Problems in decentralized decision making and computation.
Ph.D. Dissertation, MIT, 1984.

\bibitem{Tsitsiklis2}
Tsitsiklis J.~N., Bertsekas D., Athans M.
Distributed asynchronous deterministic and stochastic gradient optimization algorithms.
\emph{IEEE Transactions on Automatic Control}, 1986, \textbf{31}(9): 803--812.

\bibitem{Farid}
Farid F.~O.
Criteria for invertibility of diagonally dominant matrices.
\emph{Linear Algebra and Its Applications}, 1995, \textbf{215}: 63--93.

\bibitem{Huang}
Huang T.~Z., Yang C.~S.
\emph{Special Matrices: Analysis and Applications}.
Beijing: Science Press, 2007. (in Chinese)

\bibitem{Horn-C}
Horn R.~A., Johnson C.~R.
\emph{Matrix Analysis}. (Chinese translation.)
Beijing: China Machine Press, 2005. (in Chinese)

\bibitem{Horn-E}
Horn R.~A., Johnson C.~R.
\emph{Matrix Analysis}.
Cambridge: Cambridge University Press, 1985.

\end{thebibliography}
\end{document}